\documentclass[a4paper,10pt]{amsart}

\usepackage{amsmath,enumerate,amsfonts,amssymb,amsthm}

\theoremstyle{plain}
\newtheorem*{Theo}{\bf THEOREM}
\newtheorem{Prop}{\bf PROPOSITION}
\newtheorem*{Lem}{\bf LEMMA}
\newtheorem*{Cor}{\bf COROLLARY}
\newtheorem*{Def}{\bf DEFINITION}
\newtheorem{Ex}{\bf EXAMPLE}
\theoremstyle{remark}
\newtheorem*{Rem}{\bf Remark}

\DeclareMathOperator{\Ker}{Ker} 
\DeclareMathOperator{\pf}{pf}

\begin{document}

\title[Monge-Amp\`ere equations and generalized complex geometry]
{\bfseries Monge-Amp\`ere equations and generalized complex geometry.\\
The two-dimensional case.}
\author{Bertrand Banos}
\address{Universit\'e de Bretagne Occidentale, 6 Avenue
Victor Le Gorgeu BP 809, 29 285 Brest\\}
\email{banos@univ-brest.fr}
\date{}

\begin{abstract}
We associate an integrable generalized complex structure to each
$2$-dimensional symplectic Monge-Amp\`ere equation of divergent
type and, using the Gualtieri $\overline{\partial}$ operator, we
characterize the conservation laws and the generating function of
such equation  as generalized holomorphic objects.
\end{abstract}

\maketitle

\section*{Introduction}

A general approach to the study of non-linear partial differential
equations, which goes back to Sophus Lie, is to see a $k$-order
equation on a $n$-dimensional manifold $N^n$ as a closed subset in
the manifold of $k$-jets $J^kN$. In particular, a second-order
differential equation lives in the space $J^2N$. Neverthess, as it
was noticed  by Lychagin in his seminal paper "Contact geometry
and non-linear second-order differential equations" (\cite{L}), it
is sometimes possible to decrease one dimension and to work on the
contact space $J^1N$. The idea is to define for any differential
form $\omega\in \Omega^n(J^1N)$, a second order differential
operator $\Delta_\omega:C^\infty(N)\rightarrow \Omega^n(N)$ acting
according to the rule
$$
\Delta_\omega(f)=j_1(f)^*\omega,
$$
where $j_1(f): N\rightarrow J^1N$ is the section corresponding to
the function $f$.

The differential equations of the form $\Delta_\omega=0$ are said
to be of  Monge-Amp\`ere type because of their "hessian - like"
non-linearity. Despite its very simple description, this classical
class of differential equations attends much interest due to its
appearence in different problems of geometry or mathematical
physics. We refer to the very rich book \emph{Contact geometry and
Non-linear Differential Equations} (\cite{KLR}) for a complete
exposition of the theory and for numerous examples.

A Monge-Amp\`ere equation $\Delta_\omega=0$ is said to be
symplectic if the Monge-Amp\`ere operator $\Delta_\omega$ is
invariant with respect to the Reeb vector field. In other words,
the $n$-form $\omega$ lives actually on the cotangent bundle
$T^*N$, and symplectic geometry takes place of contact geometry.
The Monge-Amp\`ere operator is then defined by
$$\Delta_\omega(f)=(df)^*\omega.$$
This partial case is in some sense quite generic because of the
beautiful result of Lychagin which says that any Monge-Amp\`ere
equation admitting a contact symmetry is equivalent (by a Legendre
transform on $J^1N$) to a symplectic one.

We are interested here in symplectic Monge-Amp\`ere equations in
two variables. These equations are written as :
\begin{equation}{\label{MA}}
A\frac{\partial^2 f}{\partial q_1^2} + 2B
\frac{\partial^2f}{\partial q_1\partial q_2}+ C\frac{\partial^2
f}{\partial q_2^2}+ D\Big(\frac{\partial^2 f}{\partial
q_1^2}\frac{\partial^2 f}{\partial q_2^2} - \big(\frac{\partial^2
f}{\partial q_1\partial q_2}\big)^2\Big)+E=0,
\end{equation}
with $A$, $B$, $C$, $D$ and $E$ smooth functions of
$(q,\frac{\partial f}{\partial q})$. These equations correspond to
$2$-form on $T^*\mathbb{R}^2$, or equivalently  to tensors on
$T^*\mathbb{R}^2$ using the correspondence
$$\omega(\cdot,\cdot)=\Omega(A\cdot,\cdot),$$
$\Omega$ being the symplectic form on $T^*N$. In the
non-degenerate case, the traceless part of this tensor $A$
defines either an almost complex structure or an almost product
structure and it is integrable if and only the corresponding
Monge-Amp\`ere equation is equivalent to the Laplace equation or
the wave equation. This elegant result of Lychagin and Roubtsov
(\cite{LR}) is quite frustrating: which kind of integrable
geometry could we define for more general Monge-Amp\`ere equations
?

It has been noticed in \cite{Cr} that such a pair of forms
$(\omega,\Omega)$ defines an almost generalized complex structure,
a very rich concept defined recently by Hitchin (\cite{H1}) and
developed by Gualtieri (\cite{G1}), which interpolates between
complex and symplectic geometry.  It is easy to see that this
almost generalized complex structure is integrable for a very
large class of $2D$-Monge-Amp\`ere equations, the equations of
\emph{divergent type}. This observation is the starting point for
the approach proposed in this paper: the aim is to present these
differential equations as "generalized Laplace equations".

In the first part, we write down this correspondence between
Monge-Amp\`ere equations in two variable  and $4$-dimensional
generalized complex geometry.

In the second part we study the ${\overline{\partial}}$-operator
associated with a Monge-Amp\`ere equation of divergent type and we
show how the corresponding conservation laws and generating
functions can be seen as "holomorphic objects".

\section{Monge-Amp\`ere equations and Hitchin pairs}

In what follows $M$ is the smooth symplectic space
$T^*\mathbb{R}^2$ endowed with the canonical symplectic form
$\Omega$. Our point of view is local (in particulary we do not
make any distinction between closed and exact forms) but most of
the results presented here have a global version.

A primitive $2$-form is a differential form $\omega\in
\Omega^2(M)$ such that $\omega\wedge\Omega=0$. We denote by $\bot:
\Omega^k(M)\rightarrow\Omega^{k-2}(M)$ the operator $\theta\mapsto
\iota_{X_\Omega}(\theta)$, the bivector $X_\Omega$ being the
bivecor dual to $\Omega$.  It is straightforward to check that in
dimension $4$, a $2$-form $\omega$ is primitive if and only if
$\bot\omega=0$.

\subsection{Monge-Amp\`ere operators}

\begin{Def} Let $\omega$ be a $2$-form on $M$.
A $2$-dimensional submanifold $L$ is a generalized solution of the
equation $\Delta_\omega=0$ if it is bilagrangian with respect to
$\Omega$ and $\omega$.
\end{Def}

Note that a lagrangian submanifold of $T^*\mathbb{R}^2$ which
projects isomorphically on $\mathbb{R}^2$ is a graph of a closed
$1$-form $df:\mathbb{R}^2\rightarrow T^*\mathbb{R}^2$. A
generalized solution can be thought as a smooth patching of
classical solutions of the Monge-Amp\`ere equation
$\Delta_\omega=0$ on $\mathbb{R}^2$.

\begin{Ex}[Laplace equation]
Consider the $2D$-Laplace equation
$$
f_{q_1q_1}+f_{q_2q_2}=0.
$$
It corresponds to the form $\omega=dq_1\wedge dp_2-dq_2\wedge
dp_1$, while the symplectic form is $\Omega=dq_1\wedge dp_1+
dq_2\wedge dp_2$. Introducing the complex coordinates
$z_1=q_1+iq_2$ and $z_2=p_2+ip_1$, we get $\omega+i\Omega=
dz_1\wedge dz_2$. Generalized solution of the $2D$-Laplace
equation appear  then as the complex curves of $\mathbb{C}^2$.
\end{Ex}

The following theorem (so called Hodge-Lepage-Lychagin, see
\cite{L}) establishes the $1-1$ correspondence between
Monge-Amp\`ere operators and primitive $2$-forms:

\begin{Theo}
\begin{enumerate}[i)]
\item Any $2$-form admits the unique decomposition
$ \omega=\omega_0 + \lambda\omega,$ with $\omega_0$ primitive.
\item If two primitive forms vanish on the same lagrangian
subspaces, then there are proportional.
\end{enumerate}
\end{Theo}

\begin{Rem}
A Monge-Amp\`ere operator $\Delta_\omega$ is therefore uniquely
defined by the primitive part $\omega_0$ of $\omega$, since
$\lambda\Omega$ vanish on any lagrangian submanifold.  The
function $\lambda$ can be arbitrarily chosen.
\end{Rem}

Let $\omega=\omega_0+\lambda\Omega$ be a $2$-form.  We define the
tensor $A$ by $\omega=\Omega(A\cdot,\cdot)$. One has
$A=A_{0}+\lambda Id$ and
$$
A_0^2=-\pf(\omega_0)Id,
$$
where the function $\pf(\omega_0)$ is the pfaffian of $\omega_0$
defined by
$$
\omega_0\wedge\omega_0= \pf(\omega_0)\Omega\wedge\Omega.
$$
Therefore,
$$
A^2=2\lambda A -(\lambda^2+\pf(\omega_0))Id.
$$
The equation $\Delta_\omega=0$ is said to be elliptic if
$\pf(\omega_0)>0$, hyperbolic if $\pf(\omega_0)<0$, parabolic if
$\pf(\omega_0)=0$. In the elliptic/hyperbolic case, one can define
the tensor
$$
J_0=\frac{A_{0}}{\sqrt{|\pf(\omega_0)|}}
$$
which is either an almost complex structure or an almost product
structure.

\begin{Theo}[Lychagin - Roubtsov \cite{LR}]  The
following
assertions are equivalent
\begin{enumerate}[i)]
\item The tensor  $J_{0}$ is integrable.
\item The form $\omega_0/\sqrt{|\pf(\omega_0)|}$ is closed.
\item The Monge-Amp\`ere equation $\Delta_{\omega}=0$ is
equivalent (with respect to the action of local
symplectomorphisms) to the (elliptic) Laplace equation
$f_{q_1q_1}+f_{q_2q_2}=0$ or the (hyperbolic) wave equation
$f_{q_1q_1}-f_{q_2q_2}=0$.
\end{enumerate}
\end{Theo}

Let us introduce now the Euler operator and the notion of
Monge-Amp\`ere equation of divergent type (see \cite{L}).
\begin{Def}
The Euler operator is the second order differential operator
$\mathcal{E}: \Omega^2(M)\rightarrow \Omega^2(M)$ defined by
$$
\mathcal{E}(\omega)=d \bot d\omega.
$$
A Monge-Amp\`ere equation $\Delta_\omega=0$ is said to be of
divergent type if $\mathcal{E}(\omega)=0$.
\end{Def}

\begin{Ex}[Born-Infeld Equation]
The Born-Infeld equation is
$$
(1-f_t)^2 f_{xx}+2f_tf_xf_{tx} - (1+f_x^2)f_{tt}=0.
$$
The corresponding primitive form is
$$
\omega_0=(1-p_1^2)dq_1\wedge dp_2+ p_1p_2(dq_1\wedge dp_1) +
(1+p_2^2)dq_2\wedge dp_1.
$$
with $q_1=t$ and $q_2=x$. A direct computation gives
$$
d\omega_0=3(p_1dp_2 - p_2dp_1)\wedge \Omega,
$$
and then the  Born - Infeld equation is not of divergent type.
\end{Ex}

\begin{Ex}[Tricomi equation]
The Tricomi equation is
$$
v_{xx} xv_{yy}+\alpha v_x + \beta v_y + \gamma(x,y).
$$
The corresponding primitive form is
$$
\omega_0=(\alpha p _1 + \beta p_2 +\gamma(q))dq_1\wedge dq_2+
dq_1\wedge dp_2-q_2dq_2\wedge dp_1,
$$
with $x=q_1$ and $y=q_2$. Since
$$
d\omega_0=(-\alpha dq_2 + \beta dq_1)\wedge \Omega,
$$
we conclude that the Tricomi equation is of divergent type.
\end{Ex}

\begin{Lem}
A Monge-Amp\`ere equation $\Delta_{\omega}=0$ is of divergent type
if and only if it exists a function $\mu$ on $M$ such that the
form $\omega + \mu \Omega$ is closed.
\end{Lem}

\begin{proof}
Since the exterior product by $\Omega$ is an isomorphism from
$\Omega^1(M)$ to $\Omega^3(M)$, for any $2$-form $\omega$, there
exists a $1$-form $\alpha_\omega$ such that
$$
d\omega=\alpha_\omega\wedge\Omega.
$$
Since $\bot(\alpha_\omega\wedge\Omega)=\alpha_\omega$ we deduce
that $\mathcal{E}(\omega)=0$ if and only if $d\alpha_\omega=0$,
that is $d(\omega+\mu\Omega)=0$ with $d\mu=-\alpha_\omega$.
\end{proof}

Hence, if $\Delta_\omega=0$ is of divergent type, one can choose
$\omega$ being closed.   The point is that it is not primitive in
general .

\subsection{Hitchin pairs}

Let us denote by $T$ the tangent bundle of $M$ and by $T^*$ its
cotangent bundle. The natural indefinite interior product on
$T\oplus T^*$ is
$$
(X+\xi,Y+\eta)=\frac{1}{2}(\xi(Y)+\eta(X)),
$$
and the Courant bracket on sections of $T\oplus T^*$ is
$$
[X+\xi,Y+\eta]=[X,Y]+L_X\eta-L_Y\xi -\frac{1}{2}d(\iota_X\eta-
\iota_Y\xi).
$$
\begin{Def}[Hitchin \cite{H1}]
An almost generalized complex structure is a bundle map
$\mathbb{J}: T\oplus T^*\rightarrow T\oplus T^*$ satisfying
$$
\mathbb{J}^2=-1,
$$
and
$$
(\mathbb{J}\cdot,\cdot)=-(\cdot,\mathbb{J}\cdot).
$$
Such an almost generalized complex structure is  said to be
integrable if the spaces of sections of its two eigenspaces are
closed under the Courant bracket.
\end{Def}

The standard examples are
$$
\mathbb{J}_1=\begin{pmatrix} J&0\\0&-J^*\end{pmatrix}
$$
and
$$
\mathbb{J}_2=\begin{pmatrix} 0&\Omega^{-1}\\ -\Omega &
0\end{pmatrix} $$ with $J$ a  complex structure and $\Omega$ a
symplectic form.

\begin{Lem}[Crainic \cite{Cr}]
Let $\Omega$ be a symplectic form and $\omega$ any $2$-form.
Define the tensor $A$ by $\omega=\Omega(A\cdot,\cdot)$ and the
form  $\tilde{\omega}$ by
$\tilde{\omega}=-\Omega(1+A^2\cdot,\cdot)$.

The almost generalized complex structure
\begin{equation}{\label{C}}
\mathbb{J}=\begin{pmatrix} A& \Omega^{-1}\\ \tilde{\omega} & -A^*
\end{pmatrix}
\end{equation}
is integrable if and only if $\omega$ is closed. Such a pair
$(\omega,\Omega)$ with $d\omega=0$ is called a Hitchin pair
\end{Lem}

We get then immediatly the following:

\begin{Prop}
To any $2$-dimensional symplectic Monge-Amp\`ere equation of
divergent type $\Delta_\omega=0$ corresponds a Hitchin pair
$(\omega,\Omega)$ and therefore a $4$-dimensional generalized
complex structure.
\end{Prop}

\begin{Rem}
Let $L^2\subset M^4$ be a $2$-dimensional submanifold. Let
$T_L\subset T$ be its tangent bundle and $T_L^0\subset T^*$ its
annihilator. $L$ is a generalized complex submanifold (according
to the terminology of \cite{G1}) or a generalized lagrangian
submanifold (according to the terminology of \cite{BB}) if
$T_L\oplus T^0_L$ is closed under $\mathbb{J}$. When $\mathbb{J}$
is defined by $\eqref{C}$, this is equivalent to saying that $L$
is lagrangian with respect to $\Omega$ and closed under $A$, that
is, $L$ is a generalized solution of $\Delta_\omega=0$.
\end{Rem}

\subsection{Systems of first order partial differential equations}

On $2n$-dimensional manifold, a generalized complex structure
write as
$$
\mathbb{J} =\begin{pmatrix} A& \pi\\ \sigma& -A^*\\
\end{pmatrix}
$$
with different relations detailed in \cite{Cr} between the tensor
$A$, the bivector $\pi$ and the $2$-form $\sigma$. The most
oustanding being  $[\pi,\pi]=0$, that is $\pi$ is a Poisson
bivector.

In \cite{Cr}, a generalized complex structures is said to be
non-degenerate if the Poisson bivector $\pi$ is non-degenerate,
that is, if the two eigenspaces $E=\Ker(\mathbb{J}-i)$ and
$\overline{E}=\Ker(\mathbb{J}+i)$ are transverse to $T^*$. This
leads to our symplectic form $\Omega=\pi^{-1}$ and to our $2$-form
$\omega=\Omega(A\cdot,\cdot)$.

One could also take the dual point of view and study generalized
complex structure transverse to $T$. In this situation, the
eigenspace $E$ writes as
$$
E=\big\{\xi + \iota_\xi P, \xi\in T^*\otimes \mathbb{C}\big\},
$$
with $P=\pi+i\Pi$ a complex bivector. This space defines a
generalized complex structure if and only if it is a Dirac
subbundle of $(T\oplus T^*)\otimes \mathbb{C}$ and if it is
transverse to its conjugate $\overline{E}$. According to the
Maurer-Cartan type equation described in the famous paper
\emph{Manin Triple for Lie bialgebroids} (\cite{LWX}, the first
condition is
$$
[\pi+i\Pi,\pi+i\Pi]=0.
$$
The second condition says that $\Pi$ is non-degenerate.

Hence, we obtain some analog of the Crainic's result:

\begin{Def}
A Hitchin pair of bivectors is a pair consisting of two bivectors
$\pi$ and $\Pi$, $\Pi$ being non-degenerate, and satisfying
\begin{equation}{\label{bihamilt}}
\begin{cases}
[\Pi,\Pi]=[\pi,\pi]&\\
[\Pi,\pi]=0.&\\
\end{cases}
\end{equation}
\end{Def}

\begin{Prop}
There is a 1-1 correspondence between Generalized complex
structure
$$
\mathbb{J}=\begin{pmatrix} A & \pi_A\\ \sigma& -A^*\end{pmatrix}
$$
with $\sigma$ non degenerate and Hitchin pairs of bivector
$(\pi,\Pi)$. In this correspondence, we have
$$
\begin{cases}
\sigma=\Pi^{-1}\\
A=\pi\circ\Pi^{-1}\\
\pi_A= -(1+A^2)\Pi
\end{cases}
$$
\end{Prop}

\begin{Ex}
If $\pi+i\Pi$ is non-degenerate, it defines a $2$-form
$\omega+i\Omega$ which is necessarily closed (this is the complex
version of the classical result which says that a non-degenerate
Poisson bivector is actually symplectic). We find again an Hitchin
pair. So new examples occur only in the degenerate case. Note that
$\pi+i\Pi= (A+i)\Pi$, so $\det(\pi+i\Pi)=0$ if and only if $-i$ is
an eigenvalue for $A$. In dimension $4$, this implies that
$A^2=-1$ but this is not any more true in greater dimensions (see
for example the classification of pair of $2$-forms on
$6$-dimensional manifolds in \cite{LR}). Nevertheless, the case
$A^2=-1$  is interesting by itself. It corresponds to generalized
complex structure of the form
$$
\mathbb{J}=\begin{pmatrix} J&0\\
\sigma&-J^*
\end{pmatrix}
$$
with $J$ an integrable complex structure and $\sigma$ a $2$-form
satisfying $J^*\sigma= -\sigma$ and
$$
d\sigma_J=d\sigma(J\cdot,\cdot,\cdot)+
d\sigma(\cdot,J\cdot,\cdot)+d\sigma(\cdot,\cdot,J\cdot).
$$
where $\sigma_J=\sigma(J\cdot,\cdot)$ (see \cite{Cr}). Or
equivalently $\sigma+i\sigma_J$ is a $(2,0)$-form satisfying
$$
\partial(\sigma + i\sigma_J)=0.
$$
One typical example of such geometry is the so called
HyperK\"ahler geometry with torsion which is an elegant
generalization of HyperK\"ahler geometry (\cite{GP}). Unlike the
HyperK\"aler case, such geometry are always generated by
potentials (\cite{BS}).
\end{Ex}

Let us consider now an Hitchin pair of bivectors $(\pi,\Pi)$ in
dimension $4$. Since $\Pi$ is non-degenerate, it defines two
$2$-forms $\omega$ and $\Omega$, which are not necessarily closed,
and related by the tensor $A$. A generalized lagrangian surface is
a surface closed under $A$, or equivalently, bilagrangian:
$\omega|_L=\Omega|_L=0$. Locally, $L$ is defined by two functions
$u$ and $v$ satisfying a first order system
$$
\begin{cases}
a+b\frac{\partial u}{\partial x} +c\frac{\partial u}{\partial y}+
d\frac{\partial v}{\partial x}+ e\frac{\partial  v}{\partial y} +
f\det J_{u,v}\\
A+B\frac{\partial u}{\partial x} +C\frac{\partial u}{\partial y}+
D\frac{\partial v}{\partial x}+ E\frac{\partial  v}{\partial y} +
E\det J_{u,v}\\
\end{cases}
$$
with
$$
J_{u,v}=\begin{pmatrix}\frac{\partial u}{\partial x} &
\frac{\partial u}{\partial y}\\ \frac{\partial v}{\partial x}&
\frac{\partial v}{\partial y}\\
\end{pmatrix}
$$
Such a system generalizes both Monge-Amp\`ere equations and
Cauchy-Riemann systems and is called Jacobi-system (see
\cite{KLR}).

With the help of Hitchin's formalism, we understand now the
integrability  condition \eqref{bihamilt} as a "divergent type"
condition for Jacobi equations.

\section{The $\overline{\partial}$-operator}

Let us fix now a $2D$- symplectic Monge-Amp\`ere equation of
divergent type $\Delta_\omega=0$, the $2$-form
$\omega=\omega_0+\lambda\Omega$ being closed. We still denote by
$A=A_0+\lambda$ the associated tensor.

\begin{Lem}
For any $1$-form $\alpha$, the following relation holds:
\begin{equation}{\label{B}}
\alpha\wedge\omega - B^*\alpha\wedge \Omega=0
\end{equation}
with $B=\lambda-A_0$.
\end{Lem}

\begin{proof}
Let $\alpha=\iota_X\Omega$ be  a $1$-form. Since $\omega_0$ is
primitive, we get
$$
0=\iota_X(\omega_0\wedge\Omega)=(\iota_X\omega_0)\wedge\Omega+(\iota_X\Omega)\wedge
\omega_0= A_0^*\alpha\wedge\Omega+\alpha\wedge\omega_0.
$$
Therefore,
$$
\alpha\wedge\omega=\alpha\wedge\omega_0+\lambda
\alpha\wedge\Omega= (-A_0+\lambda)^*\alpha \wedge\Omega.
$$
\end{proof}

We denote by $\mathbb{J}$ the generalized complex structure
associated with the Hitchin pair $(\omega,\Omega)$. We also define
$$
\Theta=\omega-i\Omega
$$
and
$$
\Phi=\exp(\Theta)=1+\Theta+\frac{\Theta^2}{2}.
$$

\subsection{Decomposition of forms}

Using the tensor $\mathbb{J}$, Gualtieri defines a decomposition
$$
\Lambda^*(T^*)\otimes \mathbb{C}= U_{2}\oplus U_{-1}\oplus
U_{0}\oplus U_{1}\oplus U_2
$$
which generalizes the Dolbeault decomposition for a complex
structure (\cite{G1}).

Let us introduce some notations to understand this decomposition.
The space $T\oplus T^*$ acts on $\Lambda^*(T^*)$ by
$$
\rho(X+\xi)(\theta)=\iota_X\theta + \xi\wedge \theta,
$$
and this action extends to an isomorphism (the standard spin
representation) between the Clifford algebra $CL(T\oplus T^*)$ and
the space of linear endomorphisms $End(\Lambda^*(T^*))$.

\begin{Rem}
With these notations, the eigenspace $E=\Ker(\mathbb{J}-i)$ is
also defined by
$$
E=\big\{ X+\xi\in T\oplus T^*, \rho(X+\xi)(\Phi)=0\big\},
$$
\end{Rem}

\begin{Def} The space $U_k$ is defined by
$$
U_k=\rho\big(\Lambda^{2-k}\overline{E}\big)\big(\Phi\big).
$$
\end{Def}

Note that $\mathbb{J}$ identifyed with the $2$-form
$(\mathbb{J}\cdot,\cdot)$ lives in $\Lambda^2(T\oplus T^*)\subset
CL(T\oplus T^*)$. We get then an infinitesimal action of
$\mathbb{J}$ on  $\Lambda^*(T^*)$.

\begin{Lem}[Gualtieri]
$U_k$ is the $ik$-eigenspace of $\mathbb{J}$.
\end{Lem}

\begin{Rem}
We see then immediatly that $U_{-k}=\overline{U_k}$, since
$\mathbb{J}$ is a real tensor.
\end{Rem}

\begin{Prop}
\begin{enumerate}[i)]
\item $U_2=\mathbb{C}\Phi$.
\item $U_1=\big\{\alpha\wedge\Phi, \alpha\in \Lambda^1(T^*)\otimes
\mathbb{C}\big\}.$
\item $U_0=\big\{(\theta-\frac{i}{2}\bot\theta)\wedge\Phi,
\theta\in \Lambda^2(T^*)\otimes \mathbb{C}\big\}$.
\end{enumerate}
\end{Prop}

\begin{proof}
The eigenspace $\overline{E}$ is
$$
\overline{E}=\big\{X-\iota_X\overline{\Theta}, X\in T\otimes
\mathbb{C}\big\}.
$$
Now,
$$
\rho(X-\iota_X
\overline{\Theta})(\Phi)=\iota_X\Theta+\iota_X\Theta\wedge
\Theta-\iota_X\overline{\Theta}-\iota_X\overline{\Theta}\wedge\Theta=\iota_X(\Theta-\overline{\Theta})\wedge
(1+\Theta).
$$
Since $\Theta-\overline{\Theta}=-2i\Omega$ and
$X\mapsto\iota_X\Omega$ is an isomorphism between $T$ and $T^*$,
we get then the description of $U_1$.

Choose now two complex vectors $X$ and $Y$ and define
$\alpha=\iota_X\Omega$ and $\beta=\iota_Y\Omega$:
$$
\begin{aligned}
\rho\big((X-\iota_X\overline{\Theta}) & \wedge
(Y-\iota_Y\overline{\Theta})\big)\big(\Phi\big)\\
&=\rho\big(X-\iota_X\overline{\Theta}\big)\big(-2i\beta\wedge\Phi\big)\\
&=-2i\rho\big(X-\iota_X\overline{\Theta}\big)\big(\beta+\beta\wedge\Theta\big)\\
&=-2i\big(\beta(X)(1+\Theta)-\beta\wedge\iota_X\Theta-\iota_X\overline{\Theta}\wedge\beta
-\iota_X\overline{\Theta}\wedge\beta\wedge\Theta\big)\\
 &=-2i\big(\beta(X)(1+\Theta) +
 \iota_X(\Theta-\overline{\Theta})\wedge\beta\wedge (1+\Theta)-\iota_X\Theta\wedge\beta\wedge\Theta\big)\\
 &=-2i\big(\beta(X)(1+\Theta) -
 2i\alpha\wedge\beta\wedge(1+\Theta)+
 \beta\wedge\iota_X\frac{\Theta^2}{2}\big)\\
 \end{aligned}
 $$
Moreover, since $\beta\wedge\Theta^2=0$, we have
$\beta(X)\Theta^2= \beta\wedge \iota_X\Theta^2$ and then
$$
\rho\big((X-\iota_X\overline{\Theta})\wedge
(Y-\iota_Y\overline{\Theta})\big)\big(\Phi\big)=
-2i(\beta(X)-2i\alpha\wedge\beta)\wedge\Phi.
$$
But  $\bot(\alpha\wedge\beta)=-\beta(X)=\alpha(Y)$. We obtain then
the description of $U_0$.

\end{proof}

The next proposition describes the space $U_0^\mathbb{R}$ of
\emph{real} forms in $U_0$. It is a direct consequence of the
proposition above.

\begin{Prop}{\label{Ureel}}
Let $\Lambda^2_0$ be the space of (real) primitive $2$-forms. Then
$$
U_0^{\mathbb{R}}=\big\{[\theta+a(i\Omega+1)]\wedge\Phi,\text{ $
\theta\in \Lambda^2_0$ and  $a\in \mathbb{R}$}\big\}.
$$
\end{Prop}

\begin{Rem}
We have actually
$$
(\Lambda^1\oplus \Lambda^3)\otimes\mathbb{C} = U_{-1}\oplus U_{1}
$$
and
$$
(\Lambda^0\oplus \Lambda^2\oplus\Lambda^4)\otimes
\mathbb{C}=U_{-2}\oplus U_0\oplus U_2.
$$
For example, the decomposition of a $1$-form  $\alpha\in
\Lambda^1(T^*)$ is
$$
\alpha=\frac{\alpha-iB\alpha}{2}\wedge\Phi+
\frac{\alpha+iB\alpha}{2}\wedge\overline{\Phi}.
$$
\end{Rem}

This decomposition is a pointwise decomposition. Denote now by
$\mathcal{U}_k$ the space of smooth sections of the bundle $U_k$.
The Gualtieri decomposition is now
$$
\Omega^*(M)\otimes \mathbb{C}=\mathcal{U}_{-2}\oplus
\mathcal{U}_{-1}\oplus \mathcal{U}_0\oplus \mathcal{U}_1 \oplus
\mathcal{U}_2.
$$

\begin{Def}
The operator $\overline{\partial}:\mathcal{U}_k\rightarrow
\mathcal{U}_{k+1}$ is simply  $\overline{\partial}=\pi_{k+1}\circ
d$
\end{Def}

The next theorem is completely analogous to the corresponding
statement involving an almost complex structure and the Dolbeault
operator $\overline{\partial}$.

\begin{Theo}[Gualtieri \cite{G1}]
The almost generalized complex structure $\mathbb{J}$ is
integrable if and only if
$$
d=\partial + \overline{\partial}.
$$
\end{Theo}

\begin{Ex}
Let $\alpha\in \Omega^1(M)$ be a $1$-form. From
$d(\alpha\wedge\Phi)=d\alpha\wedge\Phi$ we get
$$
\begin{cases}
\overline{\partial}(\alpha\wedge\Phi)=\frac{i}{2}(\bot d\alpha)\Phi&\\
\partial(\alpha\wedge\Phi)=(d\alpha-\frac{i}{2}\bot
d\alpha)\wedge\Phi.
\end{cases}
$$
\end{Ex}

It is worth mentionning that one can also define the real
differential operator $d^\mathbb{J}=[d,\mathbb{J}]$, or
equivalently (see \cite{C})
$$
d^{\mathbb{J}}= -i(\partial-\overline{\partial}).
$$

\begin{Rem}
Cavalcanty establishes in \cite{C}, for the particular case
$\omega=0$, an isomorphism $\Xi:
\Omega^*(M)\otimes\mathbb{C}\rightarrow
\Omega^*(M)\otimes\mathbb{C}$ satisfying
$$
\Xi(d\theta)=\partial\Xi(\theta),\;\;\;
\Xi(\delta\theta)=\overline{\partial} \Xi(\theta)
$$ with  $\delta=[d,\bot]$ the symplectic codifferential.
Since $d\delta$ is the Euler operator, Monge-Amp\`ere equations of
divergent type  write as $\Delta_{\omega}=0$ with $\Xi(\omega)$
pluriharmonic on the generalized complex manifold
$\big(M^4,\exp(i\Omega)\big)$.
\end{Rem}

\subsection{Conservation laws and Generating functions}

The notion of conservation laws is a natural generalization to
partial differential equations of the notion of  first integrals.

A $1$-form $\alpha$ is a conservation law for the equation
$\Delta_\omega=0$ if the restriction of $\alpha$ to  any
generalized solution is closed. Note that conservations laws are
actually well defined up closed forms.

\begin{Ex}
Let us consider the Laplace equation and the complex structure $J$
associated with. The $2$-form  $d\alpha$ vanish on any complex
curve if and only if $[d\alpha]_{1,1}=0$, that is
$$
\overline{\partial}\alpha_{1,0} + \partial \alpha_{0,1}=0
$$
or equivalently
$$
\overline{\partial}\alpha_{1,0}= \overline{\partial}\partial \psi
$$
for some real function $\psi$. (Here $\overline{\partial}$ is the
usual Dolbeault operator defined by the integrable complex
structure $J$.) We deduce that $\alpha -d\psi = \beta_{1,0}+
\beta_{0,1}$ with $\beta_{1,0}=\alpha_{1,0}-\partial\psi$ is a
holomorphic $(1,0)$-form.

Hence, the conservation laws of the $2D$-Laplace equation are (up
exact forms)  real part of $(1,0)$-holomorphic forms.
\end{Ex}

According to the Hodge-Lepage-Lychagin theorem, $\alpha$ is a
conservation law if and only if there exist two functions $f$ and
$g$ such that $d\alpha=f\omega+ g\Omega$. The function $f$ is
called a generating function of the Monge-Amp\`ere equation
$\Delta_\omega=0$. By analogy with the Laplace equation, we will
say that the function $g$  is the conjugate function to the
generating function $f$.

\begin{Lem}
A function $f$ is a generating function if and only if $$dBdf=0.$$
\end{Lem}

\begin{proof}
$f$ is a generating function if and only if there exists a
function $g$ such that $$0=d(f\omega+g\Omega)=df\wedge\omega +
dg\wedge\Omega=(dg+Bdf)\wedge\Omega,
$$
and therefore $g$ exists if and only if $dBdf=0$.
\end{proof}

\begin{Cor}
If $f$ is a generating function  and $g$ is its conjugate  then
for any $c\in \mathbb{C}$, $L_c=(f+ig)^{-1}(c)$ is a generalized
solution of the Monge-Amp\`ere equation $\Delta_\omega=0$.
\end{Cor}

\begin{proof}
The tangent space $T_aL_c$ is generated by the hamiltonian vector
fields  $X_f$ and $X_g$. Since
$$
\Omega(BX_f,Y)=\Omega(X_f,BY)=df(BY)=Bdf(Y)=dg(Y),
$$
we deduce that $X_g=BX_f$ and therefore $L_c$ is closed under
$B=\lambda-A_0$. $L_c$ is then closed under $A_0$ and so
bilagrangian with respect to $\Omega$ and $\omega$.
\end{proof}

\begin{Ex}
A generating function of the $2D$-Laplace equation satisfies
$dJdf=0$, and hence it is the real part of a holomorphic function.
\end{Ex}

The above lemma has a nice interpretation in the Hitchin/Gualtieri
formalism:

\begin{Prop}
A function $f$ is a generating function of the Monge-Amp\`ere
equation $\Delta_\omega=0$ if and only if $f$ is a pluriharmonic
function  on the generalized complex manifold
$(M^4,\exp(\omega-i\Omega))$, that is
$$
\partial\overline{\partial} f =0.
$$
\end{Prop}

\begin{proof}
The spaces $U_1$ and $U_{-1}$ are respectively the $i$ and $-i$
eigenspaces for the infinitesimal action of $\mathbb{J}$. So
$$
\begin{aligned}
\mathbb{J}df&=\mathbb{J}\big (\frac{df - i Bdf}{2}\wedge \Phi +
\frac{df + i Bdf}{2}\wedge \overline{\Phi}\big)\\
& =i\big (\frac{df - i Bdf}{2}\wedge \Phi -
\frac{df + i Bdf}{2}\wedge \overline{\Phi}\big)\\
& = Bdf + (B^2+1)df\wedge\Omega.
\end{aligned}
$$
Moreover,
$$
d\big((B^2+1)df\wedge
\Omega\big)=d(B^2df\wedge\Omega)=d(Bdf\wedge\omega)=(dBdf)\wedge\omega.
$$
We deduce that $d\mathbb{J}df=0$ if and only if $dBdf=0$. Since
$d\mathbb{J}df= 2i\partial\overline{\partial} f$, the proposition
is proved.
\end{proof}

Decompose the function $f$ as $f=f_{-2}+f_0+f_2$. Since $\partial
f_{-2}=0$ and $\overline{\partial} f_2=0$, $f$ is pluriharmonic if
and only if $f_0$ is so. Assume that the
$\partial\overline{\partial}$-lemma holds (see \cite{C} and
\cite{G2}). Then it exists $\psi\in \mathcal{U}_{1}$ such that
$$
\overline{\partial} f_0= \overline{\partial}\partial \psi.
$$
Define then $G_0\in \mathcal{U}_0$ by $G_0=i(\partial
\psi-\overline{\partial}\overline{\psi})$. We obtain
$$
\overline{\partial}(f_0+iG_0)=0
$$
and $f_0$ appears as the real part of an "holomorphic object".
Nevertheless, this assumption is not really clear. Does the
$\partial\overline{\partial}$-lemma always hold locally ?

The following proposition gives an alternative "holomorphic
object" when the closed form $\omega$ is primitive (that is
$\lambda=0$).

\begin{Prop}
Assume that the closed form $\omega$ is primitive and consider the
real forms $U=\omega\wedge\Phi$ and  $V=(i\Omega+1)\wedge\Phi$.

A function $f$ is a generating function of the Monge-Amp\`ere
equation $\Delta_\omega=0$ with conjugate function $g$ if and only
$$
\overline{\partial}(fU-igV)=0.
$$
\end{Prop}

\begin{proof}
According to proposition \ref{Ureel}, the closed forms $U$ and $V$
live in $\mathcal{U}_0^{\mathbb{R}}$. Therefore, $d^\mathbb{J}
(fU)=-\mathbb{J}d(fU)$ and $d^\mathbb{J} (gV)=-\mathbb{J}d(gV)$.
Since $\mathbb{J}^2=-1$ on $U_{-1}\oplus U_{1}$, we get
$$
2\overline{\partial}(fU-igV)=(d-id^\mathbb{J})(fU-igV)=(1+i\mathbb{J})(dfU-d^\mathbb{J}gV).
$$
But,
$$
dfU=df\wedge \omega\wedge\Phi= df\wedge\omega,
$$
and
$$
\begin{aligned}
d^\mathbb{J}gV&=-\mathbb{J}dg\wedge V\\
& = -\mathbb{J}(idg\wedge\Omega+ dg\wedge\Phi)\\
& = -\frac{1}{2}\mathbb{J}(dg\wedge\Phi +
dg\wedge\overline{\Phi})\\
& = -\frac{i}{2}(dg\wedge\Phi - dg\wedge\overline{\Phi})\\
&= -dg\wedge\Omega.
\end{aligned}
$$
We obtain finally
$$
2\overline{\partial}(fU-igV) = df\wedge\omega + dg\wedge\Omega.
$$
\end{proof}

\begin{Ex}[Von Karman equation]
The $2D$-Von Karman equation is
$$
v_xv_{xx}-v_{yy}=0.
$$
The corresponding primitive form is
$$
\omega=p_1dq_2\wedge dp_1+dq_1\wedge dp_2,
$$
which is obviously closed. The form $U$ and $V$ are
$$
\begin{cases}
U=p_1dq_2\wedge dp_1+dq_1\wedge dp_2 + 2p_1dq_1\wedge dq_2\wedge
dp_1\wedge dp_2&\\
V=1+p_1dq_2\wedge dp_1+dq_1\wedge dp_2 + (p_1-1)dq_1\wedge
dq_2\wedge
dp_1\wedge dp_2&\\
\end{cases}
$$
\end{Ex}

\subsection{Generalized K\"ahler partners}

Gualtieri has also introduced the notion of Generalized K\"ahler
structure. This is a pair of commuting generalized complex
structure such that the symmetric product
$(\mathbb{J}_1\mathbb{J}_2)$ is definite positive. The remarkable
fact in this theory is that such a structure gives for free two
integrable complex structures and a compatible metric (see
\cite{G1}). This theory has been  used to construct explicit
examples of bihermtian structures on $4$-dimensional compact
manifolds (see \cite{H2}).

The idea is that the $+1$-eigenspace $V_+$ of
$\mathbb{J}_1\mathbb{J}_2$ is closed under $\mathbb{J}_1$ and
$\mathbb{J}_2$ and that the restriction of $(\cdot,\cdot)$ to it
is definite positive. The complex structures and the metric come
then from the natural isomorphism $V_+\rightarrow T$.

From our point of view, this approach gives us the possibility to
associate to a given partial differential equation, natural
integrable complex structures and inner products. Nevertheless, at
least  for hyperbolic equations, such inner product should have a
signature, and we have may be to a relax a little bit the
definition of generalized K\"ahler structure:

\begin{Def}
Let $\Delta_\omega=0$ be a $2D$-symplectic Monge-Amp\`ere equation
of divergent type and let $\mathbb{J}$ be the generalized complex
structure associated with. We will say that this Monge-Amp\`ere
equation admits a generalized K\"ahler partner if it exists a
generalized complex structure $\mathbb{K}$ commuting with
$\mathbb{J}$ such that the two eigenspaces of
$\mathbb{J}\mathbb{K}$ are transverse to $T$ and $T^*$.
\end{Def}

Note that a powerful tool has been done in \cite{H2} to construct
such structures:

\begin{Lem}[Hitchin] Let $\exp{\beta_1}$ and $\exp{\beta_2}$ be
two complex closed form defining generalized complex struture
$\mathbb{J}_1$ and $\mathbb{J}_2$ on $4$-dimensional manifold.
Suppose that
$$
(\beta_1-\beta_2)^2=0=(\beta_1-\overline{\beta_2})^2
$$
then $\mathbb{J}_1$ and $\mathbb{J}_2$ commute.
\end{Lem}

Let us see now on a particular case how one can use this tool.
Consider an elliptic Monge-Amp\`ere equation $\Delta_\omega=0$
with $d\omega=0$ and $\Omega\wedge\omega=0$. Assume moreover it
exists a closed $2$-form $\Theta$ such that
$$
\Omega\wedge\Theta=\omega\wedge\Theta=0
$$
and
$$
4\omega=\Omega^2+\Theta^2.
$$
Note that $\exp(\omega-i\Omega)$ and $\exp(-\omega-i\Theta)$
satisfy the conditions of the above lemma. We suppose also that
$\Theta^2=\lambda^2\Omega$ with $\lambda$ a non vanishing
function. This implies that $\omega^2=\mu^2\Omega^2$ with
$$
\mu=\frac{\sqrt{1+\lambda^2}}{2}.
$$
The triple $(\omega,\Omega,\Theta)$ defines a metric $G$ and an
almost hypercomplex structure $(I,J,K)$ such that
$$
\omega=\mu G(I\cdot,\cdot),\;\;\; \Omega=G(J\cdot,\cdot),\;\;\;
\Theta=\lambda G(K\cdot,\cdot).
$$
Define now the two almost complex structures
$$
I_+=\frac{K+\lambda J}{\mu},\;\;\;I_-=\frac{K-\lambda J}{\mu}.
$$
From
$$
\omega=\frac{\Omega+\Theta}{2}(I_-\cdot,\cdot)
$$
and
$$
\omega=\frac{\Omega-\Theta}{2}(I_+\cdot,\cdot)
$$
we deduce that $I_+$ and $I_-$ are integrable.

\begin{Lem}
A function $g$ is the conjugate of a generating function $f$ of
the Monge-Amp\`ere equation $\Delta_\omega=0$ if and only if
$$
dI_+dg=-dI_-dg.
$$
\end{Lem}

\begin{proof}
$f$ is a generating function with conjugate $g$ if and only if
$$
0=df\wedge\omega+ dg\wedge\Omega = (-\mu Kdf + dg)\wedge\Omega
$$
 that is if and only if $d\frac{K}{\mu}dg=0$.
 \end{proof}

\begin{Ex}
Consider again the Von Karman equation
$$
v_xv_{xx}-v_{yy}=0.
$$
with corresponding primitive and closed form
$$
\omega=p_1dq_2\wedge dp_1+dq_1\wedge dp_2.
$$
Define then $\Theta$ by
$$
\Theta=dp_1\wedge dp_2+(1+4p_1)dq_1\wedge dq_2.
$$
With the triple $(\omega,\Omega,\Theta)$ we construct $I_+$ and
$I_-$ defined by
$$
I_+=\frac{1}{2}\begin{pmatrix}
0&-1&1&0\\-1/p_1&0&0&-1/p_1\\
-(1+4p_1)/p_1&0&0&-1/p_1\\
0&1+4p_1&-1&0\\
\end{pmatrix}
$$
$$
I_-=\frac{1}{2}\begin{pmatrix}
0&-1&-1&0\\-1/p_1&0&0&1/p_1\\
(1+4p_1)/p_1&0&0&-1/p_1\\
0&-(1+4p_1)&-1&0\\
\end{pmatrix}
$$
It is worth mentioning that $I_+$ and $I_-$ are well defined for
all $p_1\neq 0$. But the metric $G$ is definite positive only for
$p_1<-\frac{1}{4}$.

\end{Ex}

\begin{Rem}
It would be very interesting to understand the behaviour of
generating functions and generalized solution of this kind of
Monge-Amp\`ere equations with respect to the Gualtieri metric. In
particulary, Gualtieri has introduced a scheming generalized
Laplacian $dd^*+d^*d$ (see \cite{G2}) and to know if generating
functions (which are pluriharmonic as we have seen) are actually
harmonic would give important informations on the global nature of
the solutions. This will be the object of further investigations.
\end{Rem}

\end{document}